\documentstyle[11pt,amssymb]{article}
\begin{document}

{\Large \noindent{\bf On discrete $q$-ultraspherical polynomials
and their duals}}

\bigskip

\noindent{\sc N. M. Atakishiyev${}^1$ and A. U. Klimyk${}^{1,2}$}
\medskip

\noindent ${}^1$Instituto de Matem\'aticas, UNAM, CP 62210 Cuernavaca, Morelos, M\'exico

\noindent ${}^2$Bogolyubov Institute for Theoretical Physics, 03143 Kiev, Ukraine

\medskip


\begin{abstract}
We show that a confluent case of the big $q$-Jacobi polynomials
$P_n(x;a,b,c;q):={}_3\phi_2(q^{-n},abq^{n+1},x;aq,cq;q,q)$, which
corresponds to $a=b=-c$, leads to a discrete orthogonality relation
for imaginary values of the parameter $a$ (outside of its commonly
known domain $0<a<q^{-1}$). Since $P_n(x;q^\alpha,q^\alpha,-q^\alpha;q)$
tend to Gegenbauer (or ultra\-spherical) polynomials in the limit as
$q\to 1$, this family represents yet another $q$-extension of these
classical polynomials, different from the continuous $q$-ultraspherical
polynomials of Rogers. The dual family with respect to the polynomials
$P_n(x;a,a,-a;q)$ (i.e., the dual discrete $q$-ultraspherical polynomials)
corresponds to the indeterminate moment problem, that is, these polynomials
have infinitely many orthogonality relations. We find orthogonality relations
for these polynomials, which have not been considered before. In particular,
extremal orthogonality measures for these polynomials are derived.
 \end{abstract}



\bigskip



\noindent{\bf 1. Introduction}
 \medskip

This paper deals with orthogonality relations for $q$-orthogonal
polynomials. It is well known that each family of $q$-orthogonal
polynomials corresponds to determinate or indeterminate moment
problem. If a family corresponds to determinate moment problem,
then there exists only one positive orthogonality measure $\mu$
for these polynomials and they constitute a complete orthogonal
set in the Hilbert space $L^2(\mu)$. If a family corresponds to
indeterminate moment problem, then there exists infinitely many
orthogonality measures $\mu$ for these polynomials. Moreover,
these measures are divided into two parts: extremal measures and
non-extremal measures. If a measure $\mu$ is extremal, then the
corresponding set of polynomials constitute a complete orthogonal
set in the Hilbert space $L^2(\mu)$. If a measure $\mu$ is not
extremal, then the corresponding family of polynomials is not
complete in the Hilbert space $L^2(\mu)$.

Importance of orthogonal polynomials and their orthogonality
measures stems from the fact that with each family of orthogonal
polynomials one can associate a closed symmetric (or self-adjoint)
operator $A$, representable by a Jacobi matrix. If the
corresponding moment problem is indeterminate, then the operator
$A$ is not self-adjoint and has infinitely many self-adjoint
extensions. If the operator $A$ has a physical meaning, then these
self-adjoint extensions are especially important. These extensions
correspond to extremal orthogonality measures for the same set of
polynomials and can be constructed by means of these measures
(see, for example, [1], Chapter VII). If the family of polynomials
corresponds to determinate moment problem, then the corresponding
operator $A$ is self-adjoint and its spectrum is determined by an
orthogonality relation for the polynomials. Moreover, the spectral
measure for the operator $A$ is constructed by means of the
orthogonality measure for the corresponding polynomials (see [1],
Chapter VII). Thus, orthogonality measures for polynomials are of
great interest for applications in operator theory and quantum
mechanics.

In the present paper we deal in fact with orthogonality measures
for big $q$-Jacobi polynomials and their duals. It is well known
that the big $q$-Jacobi polynomials $P_n(x;a,b,c;q)$ are
orthogonal for values of the parameters in the intervals
$0<a,b<q^{-1}$, $c<0$, and the corresponding moment problem is
determinate. We show that these polynomials are also orthogonal
outside of these intervals in a confluent case when $a=b=-c$.
Since the polynomials $P_n(x;q^\alpha,q^\alpha,-q^\alpha|q)$ tend
to ultraspherical polynomials when $q\to 1$, we call them discrete
$q$-ultraspherical polynomials (because their orthogonality
measure is discrete, contrary to the orthogonality measure for
continuous $q$-ultraspherical polynomials of Rogers). They
correspond to determinate moment problem. We give explicitly an
orthogonality relation for $P_n(x; a,a,-a|q)$ when $a$ becomes
imaginary.

In [2] we introduced a family of $q$-orthogonal polynomials
$D_n(\mu(x);a,b,c|q)$, dual to big $q$-Jacobi polynomials (they
correspond to indeterminate moment problem), and found an
orthogonality relation for them. The corresponding orthogonality
measure is not extremal. In the present paper we consider the dual
$q$-Jacobi polynomials $D_n(\mu(x);a,b,c|q)$ for $a=b=-c$. We find
two new orthogonality measures for this case when $0<a<q^{-1}$,
which are extremal. Besides, we derive two orthogonality measures
for $D_n(\mu(x);a,a,-a|q)$, when $a$ is an imaginary number. These
measures are also extremal. We also found infinitely many
orthogonality relations for these polynomials, which are difficult
to define whether they are extremal or not.

Throughout the sequel we always assume that $q$ is a fixed
positive number such that $q<1$. We use (without additional
explanation) notation of the theory of special functions (see,
for example, [3] and [4]).
\medskip

\noindent{\bf 2. Big and little $q$-Jacobi polynomials}
\medskip

The big $q$-Jacobi polynomials $P_n(x ;a,b,c;q)$, introduced by G.
E. Andrews and R. Askey [5], are defined by the formula
$$
P_n(x ;a,b,c;q):= {}_3\phi_2 (q^{-n}, abq^{n+1}, x ;\; aq,cq; \;
q,q ) . \eqno (1)
$$
The discrete orthogonality relation for these polynomials is
$$
r(a,b,c) \sum_{n=0}^\infty \frac{(aq,abq/c;q)_n
q^n}{(aq/c,q;q)_n}P_m(aq^{n+1})P_{m'}(aq^{n+1})
$$  $$
 +\, r(b,a,ab/c)
\sum_{n=0}^\infty \frac{(bq,cq;q)_n
q^n}{(cq/a,q;q)_n}P_m(cq^{n+1})P_{m'}(cq^{n+1})
$$  $$
=\frac{(1-abq)(bq,abq/c,q;q)_m}{(1-abq^{2m+1})(aq,abq,cq;q)_m}
 (-ac)^m q^{m(m+3)/2}\,\delta_{mm'}\,  ,       \eqno (2)
$$
where $r(a,b,c):=(bq, cq;q)_\infty /(abq^2,c/a;q)_\infty$. This
orthogonality relation holds for $0<a,b<q^{-1}$ and $c<0$. For
these values of the parameters the big $q$-Jacobi polynomials (1)
correspond to the determinate moment problem, that is, the
orthogonality measure in (2) is unique.

We also need below an explicit form for the little $q$-Jacobi
polynomials. They are given by the formula
$$
p_n(x ;a,b|q):={}_2\phi_1 (q^{-n}, abq^{n+1};\; aq; \; q,qx )
\eqno (3)
$$
(see, for example, (7.3.1) in [3]). The discrete orthogonality
relation for the polynomials $p_m(q^n)\equiv p_m(q^n;a,b|q)$ is of
the form
$$
\sum_{n=0}^\infty
\frac{(bq;q)_n(aq)^n}{(q;q)_n}p_m(q^n)p_{m'}(q^n)
=\frac{(abq^2;q)_\infty}{(aq;q)_\infty}
\frac{(1-abq)(aq)^m\,(bq,q;q)_m} {(1-abq^{2m+1})\,(abq,aq;q)_m}\,
\delta_{mm'}\, , \eqno (4)
$$
where $0<a<q^{-1}$ and $b<q^{-1}$. For these values of the
parameters the little $q$-Jacobi polynomials also correspond to
the determinate moment problem.

In [2] we discussed two families of orthogonal polynomials, which
are dual to big and little $q$-Jacobi polynomials, respectively.
The dual big $q$-Jacobi polynomials $D_n(\mu (x); a,b,c|q)\equiv
D_n(\mu (x;ab); a,b,c|q)$ are defined as
$$
D_n(\mu (x;ab); a,b,c|q):= {}_3\phi_2 (q^{-x},abq^{x+1},q^{-n};\;
  aq, abq/c;\; q,aq^{n+1}/c )  ,  \eqno(5)
$$
where $\mu (x;\alpha):=q^{-x}+\alpha q^{x+1}$ represents a
$q$-quadratic lattice. For $0<a,b<q^{-1}$ and $c<0$ they satisfy a
discrete orthogonality relation with respect to the measure,
supported on the points $\mu(m;ab)=q^{-m}+ab\, q^{m+1}$,
$m=0,1,2,\cdots$ (see [2], formula (29)). The polynomials $D_n(\mu
(x;ab); a,b,c|q)$ with these values of the parameters correspond
to indeterminate moment problem and the orthogonality measure for
them in [2] is not extremal.

The dual little $q$-Jacobi polynomials $d_n(\mu (x); a,b|q)\equiv
d_n(\mu (x;ab); a,b|q)$ are given by the formula
$$
d_n(\mu (x;ab); a,b|q):= {}_3\phi_1(q^{-x},ab\,q^{x+1},q^{-n};\;
bq;\; q,q^n/a)   \eqno (6)
$$
and they obey the orthogonality relation
$$
\sum_{m=0}^\infty \,\frac{(1-abq^{2m+1})(abq,bq;q)_m}
{(1-abq)(aq,q;q)_m a^{-m} q^{-m^2}}d_n(\mu(m)) d_{n'}(\mu (m))  =
\frac{(abq^2;q)_\infty}{(aq;q)_\infty}
\,\frac{(q;q)_n(aq)^{-n}}{(bq;q)_n}\, \delta_{nn'}  \eqno(7)
$$
with $\mu(m)\equiv \mu(m;ab)$, which is valid for $0<a<q^{-1}$ and
$b<q^{-1}$ (see formula (13) in [2]). For these values of the
parameters $a$ and $b$, the polynomials (6) correspond to
indeterminate moment problem and the orthogonality measure in (7)
is extremal.
\medskip

\noindent{\bf 3. Discrete $q$-ultraspherical polynomials}
\medskip

For the big $q$-Jacobi polynomials $P_n(x ;a,b,c;q)$ the following
limit relation holds:
$$
\lim_{q\uparrow 1} P_n(x ;q^\alpha,q^\beta,-q^\gamma;q)
=\frac{P^{(\alpha,\beta)}_n(x)}{P^{(\alpha,\beta)}_n(1)} ,
$$
where $\gamma$ is real. Therefore, $\lim_{q\uparrow 1} P_n(x
;q^\alpha,q^\alpha,-q^\gamma;q)$ is a multiple of the Gegenbauer
(ultraspherical) polynomial $C_n^{(\alpha-1/2)}(x)$. For this
reason, we introduce the notation
$$
C_n^{(a^2)}(x;q):=P_n(x; a,a,-a;q)={}_3\phi_2
(q^{-n},a^2q^{n+1},x;\;
  aq, -aq;\;  q,q ) .    \eqno (8)
$$
It is obvious  from (8) that $C_n^{(a)}(x;q)$ is a rational
function in the parameter $a$.

>From the recurrence relation for the big $q$-Jacobi polynomials
(1) one readily verifies that the polynomials (8) satisfy the
following three-term recurence relation:
$$
x\, C_n^{(a)}(x;q)=A_n(a)\, C_{n+1}^{(a)}(x;q) +C_n(a)\,
C_{n-1}^{(a)}(x;q), \eqno (9)
$$
where $A_n(a)=(1-aq^{n+1})/(1-aq^{2n+1})$, $C_n(a)=1-A_n(a)$ and
$C_0^{(a)}(x;q)\equiv 1$.

An orthogonality relation for $C_n^{(a)}(x;q)$, which follows from
that for the big $q$-Jacobi polynomials and is considered in the
next section, holds for positive values of $a$. We shall see that
the polynomials $C_n^{(a)}(x;q)$ are orthogonal also for imaginary
values of $a$ and $x$. In order to dispense with imaginary numbers
in this case, we introduce the following notation:
$$
\tilde C_n^{(a^2)}(x;q):= (-{\rm i})^n C_n^{(-a^2)}({\rm i}x;q)
=(-{\rm i})^n P_n({\rm i}x; {\rm i}a,{\rm i}a,-{\rm i}a;q),
 \eqno (10)
$$
where $x$ is real and $0<a<\infty$. The polynomials $\tilde
C_n^{(a^2)}(x;q)$ satisfy the recurrence relation
$$
x\tilde C_n^{(a)}(x;q)=\tilde A_n (a)\, \tilde C_{n+1}^{(a)}(x;q)
+\tilde C_n(a)\, \tilde C_{n-1}^{(a)}(x;q),  \eqno (11)
$$
where $\tilde A_n(a)=A_n(-a)=(1+aq^{n+1})/(1+aq^{2n+1})$, $\tilde
C_n(a)=\tilde A_n(a)-1$, and $\tilde C_0^{(a)}(x;q)\equiv 1$. Note
that $\tilde A_n(a) \ge 1$ and, therefore, coefficients in the
recurrence relation for $\tilde C_n^{(a)}(x;q)$ satisfy the
conditions $\tilde A_n(a)\tilde C_{n+1}(a)>0$ of Favard's
characterization theorem for $n=0,1,2,\cdots$ (see, for example,
[3]). This means that these polynomials are orthogonal with
respect to a positive measure with infinitely many points of
support. An explicit form of this measure is derived in the next
section.

So, we have
$$
\tilde C_n^{(a)}(x;q)=(-{\rm i})^n C_n^{(-a)}({\rm i}x;q)= (-{\rm
i})^n {}_3\phi_2\left( \left. {q^{-n},-aq^{n+1},{\rm i}x
  \atop {\rm i}\sqrt{a}q, -{\rm i}\sqrt{a}q} \right| q,q \right) .
  \eqno (12)
$$
{}From the recurrence relation (11) it follows that the
polynomials (12) are real for $x\in {\Bbb R}$ and $0<a<\infty$.
From (12) it is also obvious that they are rational functions in
the parameter $a$.

We show below that the polynomials $C_n^{(a)}(x;q)$ and $\tilde
C_n^{(a)}(x;q)$, interrelated by (12), are orthogonal with respect
to discrete measures. For this reason, they may be regarded as a
discrete version of $q$-ultraspherical polynomials of Rogers (see,
for example, [6]).
\medskip

{\bf Proposition 1.} {\it The following expressions for the
discrete $q$-ultraspherical polynomials (8) hold:
$$
C_{2k}^{(a)}(x;q)=\frac{(q;q^2)_k\, a^{k}}{(aq^2;q^2)_k}
(-1)^kq^{k(k+1)}p_k(x^2/aq^2; q^{-1},a|q^2),  \eqno (13)
$$   $$
C_{2k+1}^{(a)}(x;q)=\frac{(q^3;q^2)_k\, a^{k}}{(aq^2;q^2)_k}
(-1)^kq^{k(k+1)}x\, p_k(x^2/aq^2; q,a|q^2),  \eqno (14)
$$
where $p_k(y;a,b|q)$ are the little $q$-Jacobi polynomials (3).}
\medskip

{\it Proof.} To start with (13), apply Singh's quadratic
transformation for a terminating ${}_3\phi_2$ series
$$
{}_3\phi_2\left( \left. {a^2,\ b^2,\ c
  \atop abq^{1/2}, -abq^{1/2}} \right| q,q \right) =
  {}_3\phi_2\left( \left. {a^{2},\ b^2,\ c^{2}
  \atop a^2b^2q,\ 0 } \right| q^2,q^2 \right) , \eqno (15)
$$
which is valid when both sides in (15) terminate (see [3], formula
(3.10.13)), to the expression in (8) for the $q$-ultraspherical
polynomials $C_{2k}^{(a)}(x;q)$. This results in the following:
$$
C_{2k}^{(a)}(x;q)=  {}_3\phi_2\left( q^{-2k},\ aq^{2k+1},\ x^{2};\
aq^2,\ 0 ;\ q^2,q^2 \right) .
$$
Now apply to the above basic hypergeometric series ${}_3\phi_2$
the transformation formula
$$
{}_2\phi_1\left( \left. {q^{-n},\ b
  \atop c } \right| q,z \right)=\frac{(c/b;q)_n}{(c;q)_n}\,
{}_3\phi_2\left( \left. {q^{-n},\ b,\  bzq^{-n}/c
  \atop bq^{1-n}/c,\ 0 } \right| q, q \right) \eqno (16)
$$
(see formula (III.7) from Appendix III in [3]) in order to get
$$
C_{2k}^{(a)}(x;q)=\frac{(q;q^2)_k\, a^{k}}{(aq^2;q^2)_k}
(-1)^kq^{k(k+1)} {}_2\phi_1\left( q^{-2k},\ aq^{2k+1};\ q;\
q^2,x^2/a \right) .
$$
Comparing this formula with the expression for the little
$q$-Jacobi polynomials (3) one arrives at (13).

One can now prove (14) by induction with the aid of (16) and the
recurrence relation (9). Indeed, since $C_0^{(a)}(x;q)\equiv 1$
and $A_0=1$, one obtains from (9) that $C_1^{(a)}(x;q)=x$. As the
next step use the fact that $C_2^{(a)}(x;q)={}_3\phi_2
(q^{-2},aq^3,x^2;\; aq^2,0;\; q^2,q^2)$ to evaluate from (9)
explicitly that
$$
C_{3}^{(a)}(x;q)=x\, {}_3\phi_2\left( q^{-2},\ aq^5,\ x^2;\ aq^2,\
0;\  q^2, q^2 \right) .
$$
So, let us suppose that
$$
C_{2k-1}^{(a)}(x;q)=x\, {}_3\phi_2\left( q^{-2(k-1)},\ aq^{2k+1},\
x^2;\  aq^2,\ 0;\  q^2, q^2 \right) \eqno (17)
$$
for $k=1,2,3,\cdots$, and evaluate a sum $A_{2k}^{-1}x\,
C_{2k}^{(a)}(x;q)+(1-A_{2k}^{-1})C_{2k-1}^{(a)}(x;q)$.  As follows
from the recurrence relation (9), this sum should be equal to
$C_{2k+1}^{(a)}(x;q)$. This is the case because it is equal to
$$
x\, \left\{ A_{2k}^{-1}\, {}_3\phi_2\left( \left. {q^{-2k},\
aq^{2k+1},\ x^2 \atop aq^2,\ 0 } \right| q^2, q^2 \right)
 +(1-A_{2k}^{-1}){}_3\phi_2\left( \left. {q^{-2(k-1)},\
aq^{2k+1},\ x^2 \atop aq^2,\ 0 } \right| q^2, q^2 \right) \right\}
$$  $$
=x\, {}_3\phi_2\left( \left. {q^{-2k},\ aq^{2k+3},\ x^2 \atop
aq^2,\ 0 } \right| q^2, q^2 \right), \eqno (18)
$$
if one takes into account readily verified identities
$$
A_{2k}^{-1}(q^{-2k};q^2)_m+(1-A_{2k}^{-1})(q^{-2(k-1)};q^2)_m
=\frac{1-aq^{2(k+m)+1}}{1-aq^{2k+1}} (q^{-2k};q^2)_m ,
$$  $$
\frac{1-aq^{2(k+m)+1}}{1-aq^{2k+1}} (aq^{2k+1};q^2)_m=
(aq^{2k+3};q^2)_m ,
$$
for $m=0,1,2,\cdots ,k$. The right side of (18) does coincide with
$C_{2k+1}^{(a)}(x;q)$, defined by the same expression (17) with
$k\to k+1$. Thus, it remains only to apply the same transformation
formula (16) in order to arrive at (14). Proposition is proved.
 \medskip

{\it Remark.} Observe that en route to proving formula (14), we
established a quadratic transformation
$$
{}_3\phi_2\left( \left. {q^{-2k-1},aq^{2k+2},x
  \atop \sqrt{a}q, -\sqrt{a}q} \right| q,q \right)=
x\,{}_3\phi_2\left( \left. {q^{-2k},aq^{2k+3},x^2
  \atop aq^2,\ 0}  \right| q^2,q^{2} \right)   \eqno (19)
$$
for the terminating basic hypergeometric polynomials ${}_3\phi_2$
with $k=0,1,2,\cdots$. The left side in (19) defines the
polynomials $C_{2k+1}^{(a)}(x;q)$ by (8), whereas the right side
follows from the expression (18) for the same polynomials. The
formula (19) represents an extension of Singh's quadratic
transformation (15) to the case when $a^2=q^{-2k-1}$ and,
therefore, the left side in (15) terminates, but the right side
does not.
\medskip

It follows from (12)--(14) that
$$
\tilde C_{2k}^{(a)}(x;q)=\frac{(q;q^2)_k\, a^{k}}{(-aq^2;q^2)_k}
(-1)^kq^{k(k+1)}p_k(x^2/aq^2; q^{-1},-a|q^2),  \eqno (20)
$$   $$
\tilde C_{2k+1}^{(a)}(x;q)=\frac{(q^3;q^2)_k\,
a^{k}}{(-aq^2;q^2)_k} (-1)^kq^{k(k+1)}x\, p_k(x^2/aq^2; q,-a|q^2).
\eqno (21)
$$
In particular, it is clear from these formulas that the
polynomials $\tilde C_{n}^{(a)}(x;q)$ are real-valued for $x\in
{\Bbb R}$ and $a>0$.
\medskip

\noindent{\bf 4. Orthogonality relations for discrete
$q$-ultraspherical polynomials}
\medskip

Since the polynomials $C_n^{(a)}(x;q)$ are a particular case of
the big $q$-Jacobi polynomials (as it is obvious from (8)), an
orthogonality relation for them follows from (2). Setting
$a=b=-c$, $a>0$, into (2) and considering the case when $m=2k$ and
$m'=2k'$, one verifies that two sums on the left of (2) coincide
(since $ab/c=-a=c$) and we obtain the following orthogonality
relation for $C_{2k}^{(a)}(x;q)$:
$$
\sum_{s=0}^\infty \frac{(aq^2;q^2)_s q^s}{(q^2;q^2)_s}
C_{2k}^{(a)}(\sqrt{a}q^{s+1};q)C_{2k'}^{(a)}(\sqrt{a}q^{s+1};q) =
\frac{(aq^3;q^2)_\infty} {(q;q^2)_\infty}\frac{(1{-}aq)
a^{2k}}{1{-}aq^{4k+1}} \frac{(q;q)_{2k}q^{k(2k+3)}}{(aq;q)_{2k}}
\delta_{kk'},\eqno (22)
$$
where $\sqrt{a}$, $a>0$ denotes a positive value of the root.
Thus, {\it the family of polynomials $C_{2k}^{(a)}(x;q)$,
$k=0,1,2,\cdots$, with $0<a<q^{-2}$, is orthogonal on the set of
points $\sqrt{a}q^{s+1}$, $s=0,1,2,\cdots$.}

As we know, the polynomials $C_{2k}^{(a)}(x;q)$ are functions in
$x^2$, that is, $C_{2k}^{(a)}(\sqrt{a}q^{s+1};q)$ is in fact a
function in $aq^{2s+2}$. From (22) it follows that the set of
functions $C_{2k}^{(a)}(x;q)$, $k=0,1,2,\cdots$, constitute a
complete basis in the Hilbert space ${\frak l}^2$ of functions
$f(x^2)$ with the scalar product
$$
(f_1,f_2)=\sum_{s=0}^\infty \frac{(aq^2;q^2)_s\, q^s}{(q^2;q^2)_s}
f_1(aq^{2s+2})\overline{f_2(aq^{2s+2})} .
$$
This result can be also obtained from the orthogonality relation
for the little $q$-Jacobi polynomials, if one takes into account
formula (13).

Putting $a=b=-c$, $a>0$, into (2) and considering the case when
$m=2k+1$ and $m'=2k'+1$, one verifies that two sums on the left of
(2) again coincide and we obtain the following orthogonality
relation for $C_{2k+1}^{(a)}(x;q)$:
$$
\sum_{s=0}^\infty \frac{(aq^2;q^2)_s\, q^s}{(q^2;q^2)_s}
C_{2k+1}^{(a)}(\sqrt{a}\, q^{s+1};q)C_{2k'+1}^{(a)}(\sqrt{a}\,
q^{s+1};q)
 $$  $$
 =\frac{(aq^3;q^2)_\infty}
{(q;q^2)_\infty}\frac{(1-aq)\, a^{2k+1}}{(1-aq^{4k+3})}
\frac{(q;q)_{2k+1}}{(aq;q)_{2k+1}}\,  q^{(k+2)(2k+1)}\delta_{kk'}.
\eqno (23)
$$
{\it The polynomials $C_{2k+1}^{(a)}(x;q)$, $k=0,1,2,\cdots$, with
$0<a<q^{-2}$, are thus orthogonal on the set of points
$\sqrt{a}\,q^{s+1}$, $s=0,1,2,\cdots$.}

The polynomials $x^{-1}C_{2k+1}^{(a)}(x;q)$ are functions in
$x^2$, that is, $x^{-1}C_{2k+1}^{(a)}(\sqrt{a}\,q^{s+1};q)$ are in
fact functions in $aq^{2s+2}$. It follows from (23) that the
collection of functions $C_{2k+1}^{(a)}(x;q)$, $k=0,1,2,\cdots$,
constitute a complete basis in the Hilbert space ${\frak l}^2$ of
functions of the form $F(x)=xf(x^2)$ with the scalar product
$$
(F_1,F_2)=\sum_{s=0}^\infty \frac{(aq^2;q^2)_s\,
q^s}{(q^2;q^2)_s}F_1(\sqrt{a}\, q^{s+1})\overline{F_2(\sqrt{a}\,
q^{s+1})} .
$$
Again, this result can be obtained also from the orthogonality
relation for the little $q$-Jacobi polynomials if one takes into
account formula (14).

We have shown that the polynomials $C_{2k}^{(a)}(x;q)$,
$k=0,1,2,\cdots$, as well as the polynomials
$C_{2k+1}^{(a)}(x;q)$, $k=0,1,2,\cdots$, are orthogonal on the set
of points $\sqrt{a}\,q^{s+1}$, $s=0,1,2,\cdots$. However, the
polynomials $C_{2k}^{(a)}(x;q)$, $k=0,1,2,\cdots$, are not
orthogonal to the polynomials $C_{2k+1}^{(a)}(x;q)$,
$k=0,1,2,\cdots$, on this set of points. In order to prove that
the polynomials $C_{2k}^{(a)}(x;q)$, $k=0,1,2,\cdots$, are
orthogonal to the polynomials $C_{2k+1}^{(a)}(x;q)$,
$k=0,1,2,\cdots$, one has to take the set of points $\pm
\sqrt{a}\, q^{s+1}$, $s=0,1,2,\cdots$. Since the polynomials from
the first set are even and the polynomials from the second set are
odd, for each $k,k'\in \{ 0,1,2,\cdots\}$ the infinite sum
$$
I_1\equiv \sum_{s=0}^\infty \frac{(aq^2;q^2)_s\, q^s}{(q^2;q^2)_s}
\,C_{2k}^{(a)}(\sqrt{a}\,q^{s+1};q)\,C_{2k'+1}^{(a)}(\sqrt{a}\,q^{s+1};q)
$$
coincides with the following one
$$
I_2\equiv - \sum_{s=0}^\infty \frac{(aq^2;q^2)_s\,q^s}{(q^2;q^2)_s}
\,C_{2k}^{(a)}(-\sqrt{a}\,q^{s+1};q)\,C_{2k'+1}^{(a)}(-\sqrt{a}\,q^{s+1};q).
$$
Therefore, $I_1-I_2=0$. This gives the orthogonality of
polynomials from the set $C_{2k}^{(a)}(x;q)$, $k=0,1,2,\cdots$,
with the polynomials from the set $C_{2k+1}^{(a)}(x;q)$,
$k=0,1,2,\cdots$. The orthogonality relation for the whole set of
polynomials $C_{n}^{(a)}(x;q)$, $n=0,1,2,\cdots $, can be written
in the form
$$
\sum_{s=0}^\infty \sum_{\varepsilon=\pm 1} \frac{(aq^2;q^2)_s\,
q^s}{(q^2;q^2)_s}\,C_{n}^{(a)}(\varepsilon\sqrt{a}\,q^{s+1};q)\,C_{n'}^{(a)}
(\varepsilon\sqrt{a}\,q^{s+1};q)
 $$   $$
=\frac{(aq^3;q^2)_\infty} {(q;q^2)_\infty}\frac{(1-aq)\,
a^{n}}{(1-aq^{2n+1})} \frac{(q;q)_{n}\,q^{n(n+3)/2}}{(aq;q)_{n}}\,
\delta_{nn'}. \eqno (24)
$$
We thus conclude that {\it the polynomials $C_{n}^{(a)}(x;q)$,
$n=0,1,2,\cdots $, with $0<a<q^{-2}$ are orthogonal on the set of
points $\pm \sqrt{a}\,q^{s+1}$, $s=0,1,2,\cdots$.}

An orthogonality relation for the polynomials $\tilde
C_{n}^{(a)}(x;q)$, $n=0,1,2,\cdots$, is derived by using the
relations (20) and (21), as well as the orthogonality relation for
the little $q$-Jacobi polynomials. Writing down the orthogonality
relation (4) for the polynomials $p_k(x^2/aq^2; q^{-1},-a|q^2)$
and using the relation (20), one finds an orthogonality relation
for the set of polynomials $\tilde C_{2k}^{(a)}(x;q)$,
$k=0,1,2,\cdots$. It has the form
$$
\sum_{s=0}^\infty \frac{(-aq^2;q^2)_s\, q^s}{(q^2;q^2)_s}
\, \tilde C_{2k}^{(a)}(\sqrt{a}\,q^{s+1};q)\,\tilde C_{2k'}^{(a)}
(\sqrt{a}\,q^{s+1};q) $$  $$
=\frac{(-aq^3;q^2)_\infty} {(q;q^2)_\infty}\frac{(1+aq)\,a^{2k}}
{(1+aq^{4k+1})}\frac{(q;q)_{2k}}{(-aq;q)_{2k}}\, q^{k(2k+3)}\delta_{kk'}. \eqno(25)
$$
Consequently, {\it the family of polynomials $\tilde C_{2k}^{(a)}(x;q)$,
$k=0,1,2,\cdots$, is orthogonal on the set of points $\sqrt{a}\, q^{s+1}$,
$s=0,1,2,\cdots$, for $a>0$.}

As in the case of polynomials $C_{2k}^{(a)}(x;q)$, $k=0,1,2,\cdots$,
the set $\tilde C_{2k}^{(a)}(x;q)$, $k=0,1,2,\cdots$, is complete in
the Hilbert space of functions $f(x^2)$ with the corresponding scalar
product.

Similarly, using formula (21) and the orthogonality relation for
the little $q$-Jacobi polynomials $p_k(x^2/aq^2; q,-a|q^2)$, we
find an orthogonality relation
$$
\sum_{s=0}^\infty \frac{(-aq^2;q^2)_s\, q^s}{(q^2;q^2)_s}\,
\tilde C_{2k+1}^{(a)}(\sqrt{a}\, q^{s+1};q)\,\tilde C_{2k'+1}^{(a)}
(\sqrt{a}\, q^{s+1};q)$$  $$
=\frac{(-aq^3;q^2)_\infty} {(q;q^2)_\infty}\frac{(1+aq)\,
a^{2k+1}}{(1+aq^{4k+3})}\frac{(q;q)_{2k+1}}{(-aq;q)_{2k+1}}\,
q^{(k+2)(2k+1)}\delta_{kk'} \eqno (26)
$$
for the set of polynomials $\tilde C_{2k+1}^{(a)}(x;q)$,
$k=0,1,2,\cdots$. We see from this relation that {\it for $a>0$
the polynomials $\tilde C_{2k+1}^{(a)}(x;q)$, $k=0,1,2,\cdots$,
are orthogonal on the same set of points $\sqrt{a}\,q^{s+1}$,
$s=0,1,2,\cdots$.}

Thus, the polynomials $\tilde C_{2k}^{(a)}(x;q)$,$k=0,1,2,\cdots$,
as well as the polynomials $\tilde C_{2k+1}^{(a)}(x;q)$, $k=0,1,2,\cdots$,
are orthogonal on the set of points $\sqrt{a}\,q^{s+1}$, $s=0,1,2,\cdots$.
However, the polynomials $\tilde C_{2k}^{(a)}(x;q)$, $k=0,1,2,\cdots$,
are not orthogonal to the polynomials $\tilde C_{2k+1}^{(a)}(x;q)$,
$k=0,1,2,\cdots$, on this set of points. As in the previous case,
in order to prove that the polynomials $\tilde C_{2k}^{(a)}(x;q)$,
$k=0,1,2,\cdots$, are orthogonal to the polynomials $\tilde
C_{2k+1}^{(a)}(x;q)$, $k=0,1,2,\cdots$, one has to consider them
on the set of points $\pm \sqrt{a}\, q^{s+1}$, $s=0,1,2,\cdots$.
Since the polynomials from the first set are even and the
polynomials from the second set are odd, then the infinite sum
$$
I_1 \equiv \sum_{s=0}^\infty \frac{(-aq^2;q^2)_s\, q^s}{(q^2;q^2)_s}
\,\tilde C_{2k}^{(a)}(\sqrt{a}\,q^{s+1};q)\,\tilde C_{2k'+1}^{(a)}
(\sqrt{a}\,q^{s+1};q)
$$
coincides with the sum
$$
I_2\equiv - \sum_{s=0}^\infty \frac{(-aq^2;q^2)_s\,q^s}{(q^2;q^2)_s}
\,\tilde C_{2k}^{(a)}(-\sqrt{a}\,q^{s+1};q)
\,\tilde C_{2k'+1}^{(a)}(-\sqrt{a}\,q^{s+1};q).
$$
Consequently, $I_1-I_2=0$. This gives the mutual orthogonality of
the polynomials $\tilde C_{2k}^{(a)}(x;q)$, $k=0,1,2,\cdots$, to
the polynomials $\tilde C_{2k+1}^{(a)}(x;q)$, $k=0,1,2,\cdots$.
Thus, the orthogonality relation for the whole set of polynomials
$\tilde C_{n}^{(a)}(x;q)$, $n=0,1,2,\cdots $, can be written in
the form
$$
\sum_{s=0}^\infty \sum_{\varepsilon=\pm 1} \frac{(-aq^2;q^2)_s\,q^s}
{(q^2;q^2)_s}\, \tilde C_{n}^{(a)}(\varepsilon\sqrt{a}\,q^{s+1};q)
\,\tilde C_{n'}^{(a)}(\varepsilon\sqrt{a}\,q^{s+1};q)
 $$   $$
=\frac{(-aq^3;q^2)_\infty}{(q;q^2)_\infty}\frac{(1+aq)\,a^{n}}
{(1+aq^{2n+1})} \frac{(q;q)_{n}}{(-aq;q)_{n}}\, q^{n(n+3)/2}\delta_{nn'}. \eqno(27)
$$
Note that the family of polynomials $\tilde C_{n}^{(a)}(x;q)$,
$n=0,1,2,\cdots $, corresponds to the determinate moment problem.
Thus, the orthogonality measure in (27) is unique. In fact,
formula (27) extends the orthogonality relation for the big
$q$-Jacobi polynomials $P_n(x;a,a,-a;q)$ to a new domain of values
of the parameter $a$.
\medskip

\noindent{\bf 5. Dual discrete $q$-ultraspherical polynomials}
\medskip

The polynomials (5) are dual to the big $q$-Jacobi polynomials (1)
(see [2]). Let us set $a=b=-c$ in the polynomials (5), as we made
before in the polynomials (1). This gives the polynomials
$$
D_n^{(a^2)}(\mu (x;a^2)|q):=D_n(\mu (x;a^2); a,a,-a|q):=\left.
{}_3\phi_2\left({q^{-x},a^2q^{x+1},q^{-n}
  \atop aq, -aq} \right| q,-q^{n+1} \right) ,   \eqno(28)
$$
where $\mu(x;a^2)=q^{-x}+a^2q^{x+1}$. We call them dual discrete
$q$-ultraspherical polynomials. They satisfy the following
three-term recurrence relation:
$$
(q^{-x}+aq^{x+1}) D_n^{(a)}(\mu (x;a)|q)= -q^{-2n-1}(1-aq^{2n+2})
D_{n+1}^{(a)}(\mu (x;a)|q)
 $$  $$
+q^{-2n-1}(1+q)D_n^{(a)}(\mu (x;a)|q)
-q^{-2n}(1-q^{2n})D_{n-1}^{(a)}(\mu (x;a)|q).
$$

For the polynomials $D_n^{(a^2)}(\mu (x;a^2)|q)$ with imaginary
$a$ we introduce the notation
$$
\tilde D_n^{(a^2)}(\mu (x;-a^2)|q):=D_n(\mu (x;-a^2); {\rm
i}a,{\rm i}a,-{\rm i}a|q):=\left.
{}_3\phi_2\left({q^{-x},-a^2q^{x+1},q^{-n}
  \atop {\rm i}aq, -{\rm i}aq} \right| q,-q^{n+1} \right) .
    \eqno(29)
$$
The polynomials $\tilde D_n^{(a)}(\mu (x;-a^2)|q)$ satisfy the
recurrence relation
$$
(q^{-x}-aq^{x+1}) \tilde D_n^{(a)}(\mu (x;-a)|q)=
-q^{-2n-1}(1+aq^{2n+2}) \tilde D_{n+1}^{(a)}(\mu (x;-a)|q)
 $$  $$
+q^{-2n-1}(1+q)\tilde D_n^{(a)}(\mu (x;-a)|q)
-q^{-2n}(1-q^{2n})\tilde D_{n-1}^{(a)}(\mu (x;-a)|q).
$$
It is obvious from this relation that the polynomials $\tilde
D_n^{(a)}(\mu (x;-a)|q)$ are real for $x\in {\Bbb R}$ and $a>0$.
For $a>0$ these polynomials satisfy the conditions of Favard's
theorem and, therefore, are orthogonal.
 \medskip

{\bf Proposition 2.} {\it The following expressions for the dual
discrete $q$-ultraspherical polynomials (28) hold:
$$
D_n^{(a)}(\mu (2k;a)|q)=d_n(\mu (k;q^{-1}a); q^{-1},a|q^2)
={}_3\phi_1\left( \left. {q^{-2k},\ aq^{2k+1},\ q^{-2n}
  \atop aq^2 } \right| q^2, q^{2n+1} \right), \eqno (30)
$$  $$
D_n^{(a)}(\mu (2k+1;a)|q)=q^nd_n(\mu (k;qa); q,a|q^2) =q^n\,
{}_3\phi_1\left( \left. {q^{-2k},\ aq^{2k+3},\  q^{-2n}
  \atop aq^2 } \right| q^2, q^{2n-1} \right), \eqno (31)
$$
where $d_n(\mu (x;bc); b,c|q)$ are the dual little $q$-Jacobi
polynomials (6).}
\medskip

{\it Proof.} Applying to the right side of (28) the formula
(III.13) from Appendix III in [3] and then Singh's quadratic
relation (15) for terminating ${}_3\phi_2$ series, after some
transformations one obtains
$$
D_n^{(a^2)}(\mu (2k;a^2)|q)=a^{-2k}q^{-k(2k+1)}
 {}_3\phi_2\left( \left.
{q^{-2k},\ a^2q^{2k+1},\ a^2q^{2n+2}
  \atop a^2q^2,\ 0 } \right| q^2, q^{2} \right).
$$
Now apply the relation (0.6.26) from [7] in order to get
$$
D_n^{(a^2)}(\mu
(2k;a^2)|q)=\frac{(q^{-2k+1};q^2)_k}{(a^2q^2;q^2)_k}
 {}_2\phi_1\left( \left.
{q^{-2k},\ a^2q^{2k+1}
  \atop q } \right| q^2, q^{2n+2} \right).
$$
Using the formula (III.8) from [3], one arrives at the expression
for $D_n^{(a^2)}(\mu (2k;a^2)|q)$ in terms of the basic
hypergeometric function from (30), coinciding with $d_n(\mu
(k;q^{-1}a^2); q^{-1},a^2|q^2)$.

The formula (31) is proved in the same way by using the relation
(19). Proposition is proved.
\medskip

For the polynomials $\tilde D_n^{(a)}(\mu (m;-a)|q)$ we have the
expressions
$$
\tilde D_n^{(a)}(\mu (2k;-a)|q)=d_n(\mu (k;-q^{-1}a);
q^{-1},-a|q^2) ={}_3\phi_1\left( \left. {q^{-2k},\ -aq^{2k+1},\
q^{-2n} \atop -aq^2 } \right| q^2, q^{2n+1} \right), \eqno (32)
$$  $$
\tilde D_n^{(a)}(\mu (2k{+}1;-a)|q)=q^nd_n(\mu (k;-qa); q,-a|q^2)
=q^n {}_3\phi_1\left( \left. {q^{-2k}, -aq^{2k+3}, q^{-2n}
\atop -aq^2 } \right| q^2, q^{2n-1} \right). \eqno (33)
$$
It is plain from the explicit formulas that the polynomials
$D_n^{(a)}(\mu (x;a)|q)$ and $\tilde D_n^{(a)}(\mu (x;a)|q)$ are
rational functions of $a$.
\medskip

\noindent{\bf 6. Orthogonality relations for dual discrete
$q$-ultraspherical polynomials}
\medskip

An example of the orthogonality relation for $D_n^{(a^2)}(\mu
(x;a^2)|q)\equiv D_n( \mu (x;a^2);a,a,-a|q)$, $0<a<q^{-1}$, has
been discussed in [2]. However, these polynomials correspond to
indeterminate moment problem and, therefore, this orthogonality
relation is not unique. Let us find another orthogonality
relations. In order to derive them we take into account the
relations (30) and (31), and the orthogonality relation (7) for
the dual little $q$-Jacobi polynomials. By means of formula (30),
we arrive at the following orthogonality relation for
$0<a<q^{-2}$:
$$
\sum_{k=0}^\infty
\frac{(1-aq^{4k+1})(aq;q)_{2k}}{(1-aq)(q;q)_{2k}}\, q^{k(2k-1)}
D_n^{(a)}(\mu (2k)|q) D_{n'}^{(a)}(\mu (2k)|q)
 =\frac{(aq^3;q^2)_\infty}{(q;q^2)_\infty}
 \frac{(q^2;q^2)_n\,q^{-n}}{(aq^2;q^2)_n}\,\delta_{nn'} ,
$$
where $\mu(2k)\equiv \mu(2k;a)$. The relation (31) leads to the
orthogonality relation, which can be written in the form
$$
\sum_{k=0}^\infty
\frac{(1-aq^{4k+3})(aq;q)_{2k+1}}{(1-aq)(q;q)_{2k+1}}\,
q^{k(2k+1)} D_n^{(a)}(\mu (2k+1)|q)\, D_{n'}^{(a)}(\mu (2k+1)|q)
 $$   $$
 =\frac{(aq^3;q^2)_\infty}{(q;q^2)_\infty}
 \frac{(q^2;q^2)_n\, q^{-n}}{(aq^2;q^2)_n} \,\delta_{nn'},
$$
where $\mu(2k+1)\equiv \mu(2k+1;a)$ and $0<a<q^{-2}$.

Thus, we have obtained two orthogonality relations for the
polynomials $ D_n^{(a)}(\mu (x;a)|q)$, $0<a<q^{-2}$, one on the
lattice $\mu (2k;a)\equiv q^{-2k}+aq^{2k+1}$, $k=0,1,2,\cdots$,
and another on the lattice $\mu (2k+1;a)\equiv
q^{-2k-1}+aq^{2k+3}$, $k=0,1,2,\cdots$. {\it The corresponding
orthogonality measures are extremal} since they are extremal for
the dual little $q$-Jacobi polynomials from formulas (30) and (31)
(see [2]).

The polynomials $\tilde D_n^{(a)}(\mu (x;-a)|q)$ also correspond
to indeterminate moment problem and, therefore, they have
infinitely many positive orthogonality measures. Some of their
orthogonality relations can be derived in the same manner as for
the polynomials $D_n^{(a)}(\mu (x)|q)$ by using the connection
(32) and (33) of these polynomials with the dual little $q$-Jacobi
polynomials (6). The relation (32) leads to the orthogonality
relation
$$
\sum_{k=0}^\infty
\frac{(1{+}aq^{4k+1})(-aq;q)_{2k}}{(1+aq)(q;q)_{2k}}\, q^{k(2k{-}1)}
\tilde D_n^{(a)}(\mu (2k)|q)\tilde D_{n'}^{(a)}(\mu (2k)|q)
 =\frac{(-aq^3;q^2)_\infty}{(q;q^2)_\infty}
 \frac{(q^2;q^2)_n\,q^{-n}}{(-aq^2;q^2)_n}\, \delta_{nn'} ,
$$
where $\mu(2k)\equiv \mu(2k;-a)$, and the relation (33) gives rise
to the orthogonality relation, which can be written in the form
$$
\sum_{k=0}^\infty
\frac{(1+aq^{4k+3})(-aq;q)_{2k+1}}{(1+aq)(q;q)_{2k+1}}\, q^{k(2k+1)}
\tilde D_n^{(a)}(\mu (2k+1)|q)\tilde D_{n'}^{(a)}(\mu (2k+1)|q)
 $$  $$
 =\frac{(-aq^3;q^2)_\infty}{(q;q^2)_\infty}
 \frac{(q^2;q^2)_n\,q^{-n}}{(-aq^2;q^2)_n} \,\delta_{nn'},
$$
where $\mu(2k+1)\equiv \mu(2k+1;-a)$. In both cases, $a$ is any
positive number.

Thus, in the case of the polynomials $\tilde D_n^{(a)}(\mu
(x;-a)|q)$ we also have two orthogonality relations. {\it The
corresponding orthogonality measures are extremal} since they are
extremal for the dual little $q$-Jacobi polynomials from formulas
(32) and (33).

Note that the extremal measures for the polynomials $
D_n^{(a)}(\mu (x)|q)$ and $\tilde D_n^{(a)}(\mu (x)|q)$, discussed
in this section, can be used for constructing self-adjoint
extensions of the closed symmetric operators, connected with the
three-term recurrence relations for these polynomials and
representable in an appropriate basis by a Jacobi matrix (details
of such construction are given in [1], Chapter VII). These
operators are representation operators for discrete series
representations of the quantum algebra $U_q({\rm su}_{1,1})$.
Moreover, the parameter $a$ for these polynomials is connected
with the number $l$, which characterizes the corresponding
representation $T_l$ of the discrete series.
\medskip

\noindent{\bf 7. Other orthogonality relations}
\medskip

The polynomials $D_n^{(a)}(\mu (x;a)|q)$ and the polynomials
$\tilde D_n^{(a)}(\mu (x;-a)|q)$ correspond to indeterminate
moment problems. For this reason, there exist infinitely many
orthogonality relations for them. Let us derive some of these
relations for the polynomials $\tilde D_n^{(a)}(\mu (x;-a)|q)$, by
using the orthogonality relations for the polynomials (5.18) in
[8], which are (up to a factor) of the form
$$
u_n(\sinh \xi; t_1,t_2|q) ={}_3\phi_1\left( \left.{q\,e^{\xi}/t_1,
\ -q\, e^{-\xi}/t_1,\ q^{-n} \atop -q^2/t_1t_2 } \right| q,
q^{n}t_1/t_2 \right)
\eqno (34) $$
and a one-parameter family of orthogonality relations for them,
characterized by a number $d$, $q\le d <1$, are given by the
formula
$$
\sum_{n=-\infty}^\infty \frac{(-t_1q^{-n}/d, t_1q^{n}d, -t_2q^{-n}/d,
t_2q^{n}d;q)_\infty}{(-t_1t_2/q;q)_\infty} \frac{d^{4n}\,q^{n(2n-1)}
(1+d^2\,q^{2n})}{(-d^2;q)_\infty (-q/d^2;q)_\infty(q;q)_\infty}
$$  $$
\times u_r\left( (d^{-1}q^{-n}-dq^n)/2; t_1,t_2 |q \right)
u_s \left((d^{-1}q^{-n}-d\,q^n)/2; t_1,t_2 |q \right)
=\frac{(q;q)_r\,(t_1/t_2)^r}{(-q^2/t_1t_2;q)_r\, q^{r}} \delta_{rs}.\eqno (35)$$
The orthogonality measure here is positive for $t_1,t_2\in {\Bbb R}$
and $t_1t_2>0$. It is not known whether these measures are extremal
or not.

In order to use the orthogonality relation (35) for the polynomials
$\tilde D_n^{(a)}(\mu (x;-a)|q)$, let us consider the
transformation formula
$$
{}_3\phi_2\left( \left. {q^{-2k},-a^2q^{2k+1},q^{-n}
  \atop {\rm i}aq, -{\rm i}aq} \right| q,-q^{n+1} \right)=
{}_3\phi_1\left( \left. {q^{-2k},-a^2q^{2k+1},q^{-2n}
  \atop -a^2q^2}  \right| q^2,q^{2n+1} \right), \eqno (36)
$$
which is true for any nonnegative integer values of $k$. This
formula is obtained by equating two expressions (29) and (32)
for the dual discrete $q$-ultraspherical polynomials $\tilde
D_n^{(a)}(\mu (2k;-a)|q)$. Observe that (36) is still valid if one
replaces the numerator parameters $q^{-2k}$ and $-a^2q^{2k+1}$ in both
sides of it by $c^{-1}q^{-2k}$ and $-ca^2\, q^{2k+1}$, $c\in {\Bbb
C}$, respectively. Indeed, the left side of (36) represents a
finite sum:
$$
{}_3\phi_2 (q^{-n},\alpha,\beta;\ \gamma,\delta;\ q,z):=\sum_{m=0}^n
\frac{(q^{-n},\alpha,\beta;q)_m}{(\gamma,\delta,q;q)_m} z^m. \eqno (37)
$$
In the case in question $\alpha=q^{-2k}$ and $\beta=-a^2q^{2k+1}$,
so the $q$-shifted factorial $(\alpha,\beta;q)_m$ in (37) is equal
to
$$
(q^{-2k},-a^2q^{2k{+}1};q)_m {=}\prod_{j=0}^{m-1}
[1{-}a^2q^{2j{+}1}{-}q^j (q^{-2k}{-}a^2q^{2k{+}1})]
{=}\prod_{j=0}^{m-1} [1{-}a^2q^{2j{+}1}{-}q^j \mu(2k;-a^2)] ,
\eqno (38)
$$
where, as before, $\mu(2k;-a^2)=q^{-2k}-a^2q^{2k+1}$. The left
side in (36) thus represents a polynomial $p_n(x)$ in the
$\mu(2k;-a^2)$ of degree $n$. In a similar manner, one easily
verifies that the right side of (36) also represents a polynomial
$p'_n(x)$ of degree $n$ in the same variable $\mu(2k;-a^2)$.
In other words, the transformation formula (36) states that the
polynomials $p_n(x)$ and $p'_n(x)$ are equal to each other on
the infinite set of distinct points $x_k=\mu(2k;-a^2)$, $k \geq n$.
Hence, they are identical.

But the point is that
$$
(c^{-1}q^{-2k},-ca^2q^{2k+1};q)_m=\prod_{j=0}^{m-1}
[1-a^2q^{2j+1}-q^j (c^{-1}q^{-2k}-ca^2q^{2k+1})]
 $$  $$
=\prod_{j=0}^{m-1} [1-a^2q^{2j+1}-q^j \mu_c(2k;-a^2)] ,
$$
where $\mu_c(2k;-a^2)=c^{-1}q^{-2k}-ca^2q^{2k+1}$. So, the
replacements $q^{-2k}\to c^{-1}q^{-2k}$ and $a^2q^{2k+1}\to
ca^2q^{2k+1}$ change only the variable, $\mu(2k;-a^2)\to
\mu_c(2k;-a^2)$, whereas all other factors in both sides of
(36) are unaltered. Thus, our statement becomes evident.

We are now in a position to establish other orthogonality
relations for the polynomials $\tilde D_n^{(a)}(\mu (x;-a)|q)$,
different from those, obtained in section 6. To achieve this, we
use the fact that the polynomials $\tilde D_n^{(a)}(\mu (x;-a)|q)$
at the points $x^{(d)}_k:= 2k + \ln(\sqrt{aq}/d)/\ln q $ are equal
to
$$
\tilde D_n^{(a)}(\mu (x^{(d)}_k;-a)|q)= {}_3\phi_2\left( \left.
{q^{-2k}d^{-1}\sqrt{aq},\ -q^{2k}d\sqrt{aq},\ q^{-n} \atop {\rm
i}\sqrt{a}\,q,\ -{\rm i}\sqrt{a}\,q } \right| q, -q^{n+1}
\right),\eqno (39)
$$
where $\mu(x^{(d)}_k;-a)= \sqrt{aq}\left( d^{-1}\,q^{-2k} - d
\,q^{2k}\right)$. From (34) and (36) it then follows that
$$
\tilde D_n^{(a)}(\mu(x^{(d)}_k;-a)|q)=u_n\left((d^{-1}\,q^{-2k}-
d\,q^{2k})/2;\sqrt{q^3/a},\sqrt{q/a}\,|q^2\right).
$$
Hence, from the orthogonality relations (35) one obtains infinite number
of orthogonality relations for the polynomials $\tilde D_n^{(a)}(\mu(x;-a)|q)$,
which are parametrized by the same $d$ as in (35). They are of the form
$$
\sum_{n=-\infty}^\infty \frac{(-t_1q^{-2n}/d,t_1q^{2n}d, -t_2q^{-2n}/d,
t_2q^{2n}d;q^2)_\infty}{(-t_1t_2/q^2;q^2)_\infty}
 \frac{d^{4n}q^{2n(2n-1)}(1+d^2q^{4n})}{(-d^2;q^2)_\infty
(-q^2/d^2;q^2)_\infty(q^2;q^2)_\infty}
 $$  $$
\times \tilde D_r^{(a)}(\mu(x^{(d)}_n;-a)|q)\tilde
D_s^{(a)}(\mu(x^{(d)}_n;-a)|q)
=\frac{(q^2;q^2)_r}{(-q^2a;q^2)^2_r}\,\delta_{rs}, \eqno (40)
$$
where $t_1=\sqrt{q^3/a}$ and $t_2=\sqrt{q/a}$.

It is important to know whether an orthogonality measure for a family
of polynomials is extremal or not. The extremality of the measures in
(40) for the polynomials $\tilde D_n^{(a)}(\mu(x;-a)|q)$ depends on the
extremality of the orthogonality measures in (35) for the polynomials
(34). If any of the measures in (35) is extremal, then the corresponding
measure in (40) is also extremal.

\medskip

\noindent{\bf 8. Concluding remarks}
\medskip

The Askey scheme [7] conveniently embraces, up to ${}_4\phi_3$-level,
all known families of orthogonal basic hypergeometric polynomials:
from continuous $q$-Hermite, Stieltjes--Wigert, and discrete $q$-Hermite
polynomials on the ground level of this scheme (for these families do
not contain any parameter other than $q$) up to the four-parameter
Askey--Wilson and $q$-Racah polynomials on the highest, fourth level.
The members of this hierarchy are known to possess the simple property:
zero values and limit cases of the parameters for any $q$-family lead
to other sets in the same hierarchy and does not yield anything novel.
The situation seems to be different for confluent cases of parameters.
So in the present paper we have studied a confluent case for the big
$q$-Jacobi polynomials $P_n(x;a,b,c;q)$ when $a=b=-c$. It turns out
that the emerging one-parameter family of $q$-polynomials represents
a $q$-extension of Gegenbauer (or ultraspherical) polynomials with a
discrete orthogonality relation with respect to the measure, supported
on the infinite set of points $\pm \, a^{1/2}\, q^{k+1}$, $k=0,1,2,\cdots$.
This family is different from all known one-parameter sets, which occupy
the first level in the Askey $q$-scheme.

The discrete $q$-ultraspherical polynomials and their duals are
evidently interesting by themselves. Besides, this particular example,
considered by us, demonstrates that other confluent cases from the
Askey $q$-scheme certainly deserve our close attention.

Another instance of similar interest is provided by the
$q$-Meixner--Pollaczek polynomials
$$
P_n(x;a|q):= e^{-{\rm i}n\phi}\,\frac{(a^2;q)_n}{a^n\,(q;q)_n}
\,\, {}_3\phi_2\left( \left. {q^{-n},a e^{{\rm i}(\theta+ 2\phi)},
a e^{-{\rm i}\theta}\atop a^2, 0} \right| q,q \right), \quad
x = \cos (\theta+\phi).$$
They are orthogonal on the interval $-\pi \leq \theta \leq \pi$ for $0<a<1$
(see, for example, section 3.9 in [7]). But let us replace $a$ by ${\rm i}a$
and assume that $\phi=-\pi/2$. One obtains then real polynomials
$$
\tilde P_n(\sin\theta;a|q):= \frac{(-a^2;q)_n}{a^n\,(q;q)_n}\,
\,{}_3\phi_2\left( \left. {q^{-n},{\rm i}a e^{-{\rm i}\theta},
 -{\rm i}a e^{{\rm i}\theta} \atop -a^2, 0} \right| q,q \right)\,,  \eqno (42)
$$
which satisfy the three-term recurrence relation
$$
2x \tilde P_n(x;a|q)=(1-q^{n+1})\tilde P_{n+1}(x;a|q)
+(1+a^2q^{n-1})\tilde P_{n-1}(x;a|q).
$$
It is now obvious that the polynomials (42) satisfy the conditions
$A_n\,C_{n+1}> 0$, $n=0,1,2,...$, of Favard's theorem for arbitrary
real $a$ and, therefore, they are orthogonal. An orthogonality
relation has the form
$$
\frac{1}{2\pi}\int_{-\pi}^\pi \tilde P_m(\sin\theta;a|q)
\tilde P_n(\sin\theta;a|q)\,w(\theta)d\theta =
(q;q)_n^{-1}(q,-a^2\,q^n;q)_\infty^{-1} \delta_{mn},
$$
where
$$
w(\theta)=| (-e^{2{\rm i}\theta};q)_\infty /
(-a^2e^{2{\rm i}\theta};q^2)_\infty |^2.
$$
A more detailed discussion of this orthogonality property will be
given elsewhere.

\medskip

\noindent{\bf Acknowledgments}
\medskip

This research has been supported in part by the SEP-CONACYT
project 41051-F and the DGAPA-UNAM project IN112300 "{\'O}ptica
Matem\'atica". A.~U.~Klimyk acknowledges the Consejo Nacional de
Ciencia y Technolog\'{\i}a (M\'exico) for a C\'atedra Patrimonial
Nivel II.

\end{document}